\documentclass{elsart}%
\usepackage{amsmath}
\usepackage{graphicx}
\usepackage{float}
\usepackage{indentfirst}
\usepackage{verbatim}
\usepackage{algorithm}
\usepackage{algpseudocode}
\usepackage{float}
\usepackage{multicol}
\usepackage{listings}
\usepackage{subcaption}
\usepackage{epsfig}
\usepackage{amsfonts}
\usepackage{amssymb}
\usepackage{epstopdf}
\usepackage{color}
\usepackage{epstopdf}%
\newtheorem{theorem}{Theorem}

\newtheorem{definition}[theorem]{Definition}

\newcommand{\splitatcommas}[1]{%
	\begingroup
	\begingroup\lccode`~=`, \lowercase{\endgroup
		\edef~{\mathchar\the\mathcode`, \penalty0
			\noexpand\hspace{0pt plus 1em}}%
		}\mathcode`,="8000 #1%
	\endgroup
}

\begin{document}
\begin{frontmatter}
\title{Multiscale DeepONet for Nonlinear Operators in Oscillatory Function Spaces for Building Seismic Wave Responses}
\author[SMU]{Lizuo Liu}
\author[SMU]{Wei Cai}

\address[SMU]{ Department of Mathematics, Southern Methodist University, Dallas, TX 75275, USA}

\bigskip
{\bf Suggested Running Head:}
\\
Multiscale DeepONet for Operators in Oscillatory Function Spaces
\\
\bigskip
{\bf Corresponding Author: }
\\
Prof. Wei Cai \\
Department of Mathematics, \\
Southern Methodist University, \\
Dallas, TX 75275\\
Email: cai@smu.edu

\begin{abstract}
In this paper, we propose a multiscale DeepONet to represent nonlinear operator between Banach spaces of highly oscillatory continuous functions. The multiscale deep neural network (DNN) utilizes a multiple scaling technique to convert high frequency function to lower frequency functions before using a DNN to learn a specific range of frequency of the function. The multi-scale concept is integrated into the DeepONet which is based on a universal approximation theory of nonlinear operators. The resulting multi-scale DeepONet is shown to be effective to represent building seismic response operator which maps oscillatory seismic excitation to the oscillatory building responses.
\end{abstract}
 
\end{frontmatter}

\textbf{}
\bigskip

\section{Introduction}

Computing operators between physical quantities defined in function spaces have many applications in forward and inverse problems in scientific and engineering computations. For example, in wave scattering in inhomogeneous or random media, the mapping between the media physical properties, which can be modelled as a random field, and the wave field is a nonlinear operator, which embodies some of most challenging computational tasks in medical imaging, geophysical and seismic problems. Another example in earthquake safety of buildings and structures, the response of structures to the seismic waves gives rise to another important operator between spaces of highly oscillatory functions. Due to the infinite dimension of function operators, active researches have been done recently in efficient representation of both forward and inverse operators with reduced models and model reductions. It is all nature that deep neural network has been considered as a machine learning tool to learn these operators.  DNNs have been used for the field of scientific and engineering
computing has shown much promise \cite{weinan18}\cite{han18}\cite{gk19}
\cite{gk20} \cite{cai19}\cite{cai20}. Recently, the DeepOnet\cite{luDeepONet2020} was constructed based on the theory of Universal Approximation of Nonlinear Operators \cite{tianpingchen1995} for learning operators between Banach spaces. Other approach for learning operators include
Fourier neural operator \cite{li2021}. In this paper,we will further develop the DeepONet for operators between highly oscillatory function spaces such as those encountered in the building seismic wave response problems and a multiscale DeepONet based on multiscale DNN concept will be investigated.

The traditional way to get the responses of buildings is by the finite element method\cite{OpenSystemEarthquake}, which could take \(O\left( 10^{2} \right) \) to \(O\left( 10^{3} \right)\) seconds for the computation of single case, depending on the grids of the finite elements. The method we developed uses the multiscale deep neural network\cite{cai20} and the DeepONet to learn the mapping between the seismic excitations and the corresponding displacements of the building at all levels excited by the seismic waves, which could give the displacements in \(O\left( 10^{-2} \right)\) seconds with relative accurate results.

The paper is constructed as follows. In section \ref{operator}, we introduce the theoretical basis of the seismic wave response operator and the idea of data augmentation of the linear differential equation case. In section \ref{multiscale}, we give a short review of the multiscale deep neural network and the integration of multiscale idea into the DeepONet structure. In section \ref{preprocess}, we introduce the data processing procedure, including the computation of responses and the elimination of the aliasing of the seismic records. The numerical results are shown in section \ref{sec:Results} and finally, a conclusion is given in section \ref{sec:conclusion}.

\section{Building seismic wave response operator $\cal R$}
\label{operator}
\noindent {\bf Building response to seismic waves} The problem we want to study is the dynamics behavior of buildings to the seismic waves\cite{NumericalSimulationLargescale2006}. We want to predict the response of the buildings instantly given the seismic excitations over a time period $[0, T]$, which defines an operator
\begin{equation}
\label{Roperator}
\mathcal{R}: P(t) \mapsto x(t) , \quad t \in [0, T]
\end{equation}
where \(P\) is the seismic wave excitation and \(x\) is the displacements of the building  excited by the seismic wave excitation \(P\).
We assume that the dynamics of a deformed body  can be described by the following equation of motion after an appropriate finite element discretization:
    \begin{equation}
         \mathbf{M} \ddot{x}\left( t \right) = P - F + H - C\dot{x}
    \label{Solid}
    \end{equation}
    with initial conditions
    \begin{equation}
    \left\{
    \begin{aligned}
        &x\left( 0 \right) = 0 \\
        &x^{\prime}\left( 0 \right) = 0
    \end{aligned}\right.,
    \label{initial}
    \end{equation}
where \(\mathbf{M}\) is the global mass matrix, \(P\) accounts for the global load vector (nodal load, body force, surface force, etc.),  \(C\) is the damping matrix, \(H\) is the global hour-glass resisting force vector handling the hour-glass deformation modes, and \(F\) is the assembly of equivalent nodal force vectors from all the elements
\begin{equation}
    F = \sum_{i=1}^{n} \int_{V_{i}}^{} \mathbf{B}^{T}\iota dV
\label{Force}
\end{equation}
where \(\mathbf{B}\) is the strain-displacement matrix, \(\iota\) is the stress vector, and \(V\) is the volume in the current configuration.

The \(\mathbf{M}, C, F\) and \(H\) are given for the operator fitting problem. The inputs and outputs of the neural network should be the global load vector \(P\left( t \right)\) and the displacements of \(x\left( t \right)\) at specific moments \(t_1,t_2,\ldots , t_{n}\).

If we only consider \(1\) level of the building and only \(1\) node in that level, \(x\) in equation \(\left( \ref{Solid} \right)\) is a 1D time series, otherwise \(x\) is an mD time series if we only set 1 node but considering \(m\) levels.

\noindent {\bf The superposition of the solutions}  Consider the training dataset \(\{ \left( p_{i}, x_{i} \right), i \in \mathcal{I}\} \), where the seismic wave excitations \(\{p_{i}\}\) and the corresponding displacements \(\{x_{i}\}\) satisfy the equation \(\left(\ref{Solid}\right)\), we build every training data to be fed to the neural network by

  \begin{equation}
    \left\{
    \begin{aligned}
        &\mathcal{P} = \sum_{i \in \iota }^{} w_{i} p_{i}, \\
        &\mathcal{X} = \sum_{i \in \iota }^{} w_{i} x_{i}
    \end{aligned}\right.
    \label{initial}
    \end{equation}
    where \(\iota\) is a randomly sampling index subset from the index set \(\mathcal{I}\) whose size is fixed and \(w_{i}\)'s are the random weights satisfying
    \begin{equation}
    \sum_{i \in  \iota}^{} w_{i} = 1
    \label{randomW}
    \end{equation}
Thus, \(\mathcal{P}\) and \(\mathcal{X}\) satisfy
   \begin{equation}
    \left\{
    \begin{aligned}
        &\mathbf{M} \ddot{\mathcal{X}} = \mathcal{P} - F + H - C\dot{\mathcal{X}},\\
        &\mathcal{X}\left( 0 \right) = 0, \\
        &\mathcal{X}^{\prime}\left( 0 \right) = 0
    \end{aligned}\right.
    \label{WeightedSol}
    \end{equation}
     by the superposition principle of the nonhomogeneous case, which construct the general idea of data augmentation for the case that there is a linear mapping relation between inputs and outputs.

\section{Multiscale DeepONet for operators in oscillatory function spaces}
\label{multiscale}
    Solutions for wave propagation and scattering in inhomogeneous media
can be viewed through an operator between the media physical properties such
as permittivity, density, etc and the wave fields through their governing
equations including elasticity wave equations, Maxwell's equations, and
acoustic wave equations. In high frequency wave problems, the operator will
have either its domain or range involves functions of highly oscillatory nature.

We will extend the DeepONet in order to represent the operator $R$ between highly oscillatory functions, namely, the seismic excitations, and the building responses, both of which are observed to include a wide range of frequencies. First of all, we define the so-called Tauber-Wiener function to be used as activation function in DNN.
    \begin{definition}
        If a function \(\sigma: \mathbb{R}\to \mathbb{R}\) \(\left( \text{continuous or discontinuous} \right)\) satisfies that all the linear combinations \(\sum_{i=1}^{I} c_{i}\sigma\left( \lambda_{i}t + \theta_{i} \right), \lambda_{i} \in \mathbb{R}, \theta_{i} \in \mathbb{R}, c_{i} \in \mathbb{R}, i=1,\ldots ,I\) are dense in every \(C\left[ a,b \right] \), then \(\sigma\) is called a Tauber-Wiener \(\left( \text{TW} \right)\) function.
    \end{definition}

The DeepONet is based on the universal approximation for nonlinear operators.
This theory gives a constructive procedure for
approximating nonlinear operator between continuous functions $u(x)$ in $x\in
K_{1}\subset\mathcal{X}$ and continuous functions in $y\in K_{2}%
\subset\mathcal{Y=}$ $\mathbb{R}^{n}$,
\begin{equation}
\mathcal{G}:u=u(x)\in V\subset C(K_{1})\rightarrow\mathcal{G(}u)(y)\in
C(K_{2})
\end{equation}
based on two uniform approximation results:

\begin{itemize}
\item \textbf{Uniform approximation of functions:} Given any $\varepsilon
_{1},$ functions $\mathcal{G(}u)(y)$ selected from a compact subset in
$C(K_{2})$ can be uniformed approximated by a two-layered neural network with
any Tauber-Wiener $\left(  \text{TW}\right)  $\ activation function
$\sigma_{t}$
\begin{equation}
|\mathcal{G(}u)(y)-\sum_{k=1}^{N}c_{k}\cdot\sigma_{t}\left(  \omega_{k}\cdot
y+\zeta_{k}\right)  |\leq\varepsilon_{1},\text{ \ }\forall\text{ }y\in K_{2}
\label{func_appr}%
\end{equation}
where $c_{k}=c_{k}(u)\circeq c_{k}(\mathcal{G(}u))$, and $\omega_{k},\zeta
_{k}$ are all independent of $y$ and the function $\mathcal{G(}u)(y)$ being approximated.
\bigskip

\item \textbf{Uniform approximation of functionals by a finite dimensional
nonlinear operator:} Given any $\varepsilon_{2},$ functional $c_{k}(u),u\in V$
(compact subset of $C(K_{1})$) can be approximated by a two-layer neural
network with any Tauber-Wiener $\left(  \text{TW}\right)  $ activation
function $\sigma_{b}$
\begin{equation}
|c_{k}(u)-\sum_{i=1}^{M}c_{i}^{k}\sigma_{b}\left(  \sum_{j=1}^{m}\xi_{ij}%
^{k}u\left(  x_{j}\right)  +\theta_{i}^{k}\right)  |\leq\varepsilon
_{2},\forall u\in V, \label{functional_appr}%
\end{equation}
where the cofficients $c_{i}^{k},\xi_{ij}^{k},\theta_{i}^{k}$ and nodes
$\{x_{j}\}_{j=1}^{m}$ and $m,M$ are all independent of $u$.
\end{itemize}

\bigskip Combining these two uniform approximation, we have an universal
approximation of nonlinear operator when restricted to a compact subset of
continuous function defined on a compact domain $K_{1}$. Namely, given any
$\varepsilon,$ we can find $\{x_{j}\}_{j=1}^{m}\subset K_{1}\subset
\mathcal{X},$coefficents $c_{i}^{k},\xi_{ij}^{k},\omega_{k},\zeta_{k}$-all
independent of continuous functions $u$ $\in V\subset C(K_{1})$ and $y\in
K_{2}\subset$ $\mathbb{R}^{n}$%
\begin{equation}
|\mathcal{G}(u)(y)-\sum_{k=1}^{N}\sum_{i=1}^{M}c_{i}^{k}\sigma_{b}\left(
\sum_{j=1}^{m}\xi_{ij}^{k}u\left(  x_{j}\right)  +\theta_{i}^{k}\right)
\cdot\sigma_{t}\left(  \omega_{k}\cdot y+\zeta_{k}\right)  |<\epsilon.
\label{operator_appro}%
\end{equation}

\bigskip

\begin{figure}[htpb]
    \begin{center}
        \includegraphics[width=0.5\textwidth]{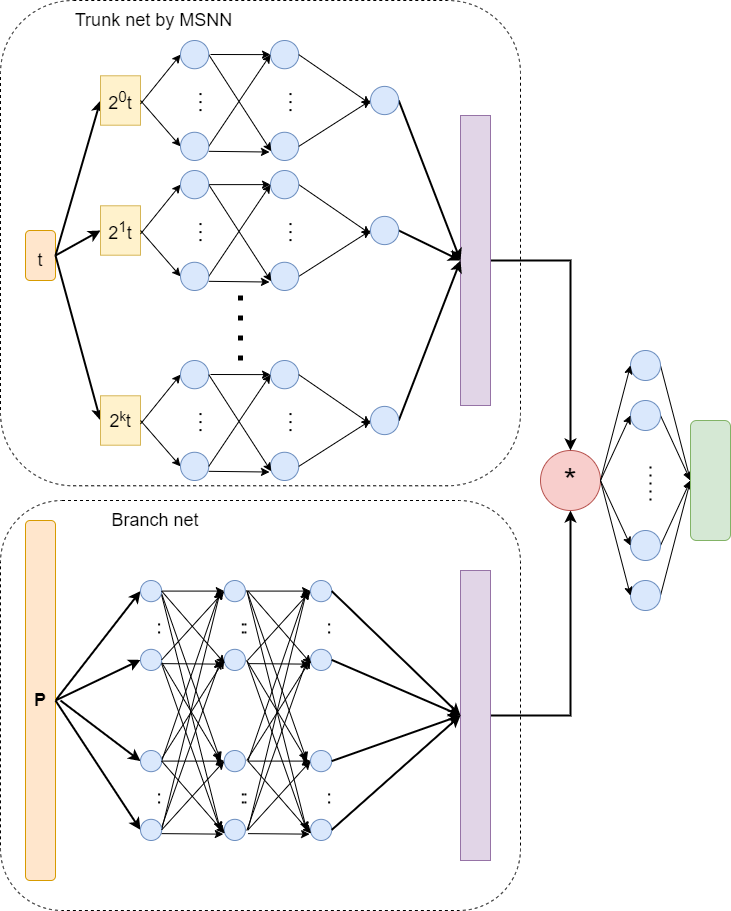}
    \end{center}
    \caption{The schematics of the (bFCN-tMS) structure.}
    \label{fig:onet_MS}
    \end{figure}

\noindent\textbf{Multiscale DeepONet} The outer sum of N-DNNs in
(\ref{operator_appro}) is identified as trunk in the DeepONet with activation
function $\sigma_{t}$, which is used for providing uniform approximation of time
dependent functions, i.e. building response in our case, sampled from a
compact subsect of $C^{0}[0,T]$; the inner sum of M DNNs (\ref{operator_appro}%
) \ with activation function $\sigma_{b}$ is used for providing uniform
approximation to functionals of infinite dimensions, selected from a compact
subspace of functionals, by m-dimensional nonlinear operators. As the building
response can be highly oscillatory, therefore, we will apply the multscale DNN
concept to generate a multiscale DeepONet as shown in Fig. \ref{fig:onet_MS}.

 In principle, we can also apply the multiscale DNN concept to the branch network, thus , we have several possible ways to combine the DeepONet and the multiscale deep neural network. The main criteria that whether we need to set specific scales for the input is whether there are high frequency components with respect to the specific variables. From the deep learning viewpoint, The DeepONet is a network with \(m+n\) inputs and only 1 output, where \(m\) is the number of the inputs of the branch net, \(n\) is the number of the inputs of the trunk net. The inputs of the trunk net \([P\left( t_1 \right), P\left( t_2 \right),\ldots P\left( t_m \right)]\) are the evaluations at the scattered sensors \(t_1,t_2,\ldots ,t_{m}\). If we want to analyse the frequency of response w.r.t every input components to determine whether the branch net should be set as multiscale scheme, we need to calculate the fourier transform with respect to each \(P\left( t_i \right)\). To do so, we need to interpolate a \(m+n\) dimension function \(f: \mathbb{R}^{m+n}\to \mathbb{R}\). If a large size of sensors are used, the analysis could be troublesome due to the difficulty of interpolating high dimension data.

 The DeepONet approximation for the seismic response operator in (\ref{Roperator}) for one location $x(t)$ in the building will be given as for $t \in [0, T]$,
 \begin{equation}
 \label{1floorR}
 x(t)=\mathcal{R}(P)(t) \approx \sum_{k=1}^{N}\sum_{i=1}^{M}c_{i}^{k}\sigma_{b}\left(
\sum_{j=1}^{m}\xi_{ij}^{k}P\left(  t_{j}\right)  +\theta_{i}^{k}\right)
\cdot\sigma_{t}\left(  \omega_{k}\cdot t+\zeta_{k}\right).
\end{equation}

 \noindent{\bf Multiscales in DeepONet} To simplify this procedure, we consider designing pre-experiments for different combinations of the DeepONet and the Multiscale DeepONet. The combinations we are interested in are
    \begin{itemize}
        \item (bMS-tFCN) The branch net is multiscale neural network, the trunk net is the general fully connected neural network.
        \item (bFCN-tMS) The branch net is fully connected neural network, the trunk net is the multiscale neural network.
        \item (bMS-tMS) The branch net and the trunk net are both multiscale neural network.
        \item (bFCN-tFCN) The branch net and the trunk net are both fully connect network.

    \end{itemize}
    The results can be seen in the section \ref{subsec:preresults}. We will select the bMS-tFCN structure based on the results of pre-experiments. The schematics of the bFCN-tMS structure is shown in Fig. \ref{fig:onet_MS}.

\noindent{\bf Loss function}
Given batch size \(\mathcal{N}\) for training process and the total number of the test records \(N\), the loss function is defined as
\begin{equation}
    \mathcal{L}\left( \theta \right) = \frac{1}{\mathcal{N}}\sum_{i=1}^{\mathcal{N}} \frac{1}{\text{max}_{j} |y_{ij}|} \Delta t \sum_{j=1}^{m} \left( f_{N N}\left( t_{j},\mathbf{u}_i,\theta \right) - y_{ij} \right)^2
\label{loss}
\end{equation}
where \(m\) is the length of signal. The penalty is set to be the reciprocal of the maximum of the absolute value of the response since the discrepancy regarding the magnitude of responses could be large corresponding to different seismic records. The larger penalty should applied to assure that the neural network could predict the response whose magnitude is smaller well.

To evaluate the training process, the mean of the relative L2 error is considered
\[
    \mathcal{L}\left( \theta \right) = \frac{1}{\mathcal{N}}\sum_{i=1}^{\mathcal{N}}\frac{ \sqrt{\Delta t \sum_{j=1}^{m} \left( f_{N N}\left( t_{j},\mathbf{u}_i,\theta \right) - y_{ij} \right)^2}}{\sqrt{\Delta t \sum_{j=1}^{m} \left( y_{ij} \right)^2}}
\]
The relative L2 error in a complete epoch is defined as
\begin{equation}
    \mathcal{L}^{\mathcal{R}}_{\text{train}} = \frac{1}{B} \sum_{k=1}^{B} \frac{1}{\mathcal{N}}\sum_{i=1}^{\mathcal{N}}\frac{ \sqrt{\Delta t \sum_{j=1}^{m} \left( f_{N N}\left( t_{j},\mathbf{u}_i,\theta^{\left( k \right)} \right) - y_{ijk} \right)^2}}{\sqrt{\Delta t \sum_{j=1}^{m} \left( y_{ijk} \right)^2}}
    \label{L_2_train}
\end{equation}
where \(B\) is the number of batches, \(\theta^{\left( k \right)}\) means the parameters of neural network at  \(k\)-th batch.

Similarly, we define the relative L2 error for the testing dataset
\begin{equation}
\mathcal{L}^{\mathcal{R}}_{\text{test}} = \frac{1}{{N}}\sum_{i=1}^{N}\frac{ \sqrt{\Delta t \sum_{j=1}^{m} \left( f_{N N}\left( t_{j},\mathbf{u}_i,\theta \right) - y_{ij} \right)^2}}{\sqrt{\Delta t \sum_{j=1}^{m} \left(  y_{ij} \right)^2}}
    \label{L_2_test}
\end{equation}
Note the \(N\) is the total number of the test records.

\section{Data preparation}
\label{preprocess}
To compute the responses of the building, we use the package openseespy\cite{OpenSystemEarthquake} and use Newmark's method to solve equation \(\left( \ref{Solid} \right)\) numerically. We need to solve \(N\) differential equations for the nodal displacements \(x\), and those \(N\) differential equations could be transformed to several main equations by expressing the displacements to the combination of the first few natural vibration modes \(\phi_{n}\). The \(\phi_{n}\)'s should be solved by the Newmark's Method. Thus we could obtain \(x\) by the inverse transform with respect to \(\phi_{n}\). The details of solving the dynamics are introduced in Chapter \(16\) in \cite{chopraDynamicsStructures2011}.

Before solving the response, we need generate simulated seismic records since the real seismic records are limited to obtain. We use the SeisMonGen\cite{WaveletbasedSimulation2017} package to generate the simulated seismic records. SeisMoGen is a software that generate earthquake ground motion simulation based on the wavelet decomposition under Priestley process assumption. For each specific real records, the generator can generate 50 random samples.

The real seismic records should be downsampled before simulating the generated samples by openseespy since the input size of the branch net is always fixed and the same time step should be considered for different seismic records. The original idea to re-sample one seismic record \(\{p_{n}\} \) is retaining every \(K\)th sample from the records such that the new records satisfies \(\tilde{f}_{n} = f_{nK}\), where \(\tilde{f}_{n}\) is the down sampled sequence. However, it is not valid to directly down-sample a sequence unless it is known that the spectrum is \(0\) at frequencies at or above the Nyquist frequency defined by the lower sampling frequency. Since this will cause the phenomenon called aliasing. Aliasing occurs when a signal is under-sampled. If the signal sampling rate is small, there will be aliasing introduced.

 \noindent\textbf{The aliasing}
    To discuss the relation of aliasing and downsampling, we follow the notations in Julius O. Simth's book\cite{MDFT}, define
    \begin{equation}
     \begin{aligned}
\operatorname{DOWNSAMPLE}_{L, m}(x) & \triangleq x(m L) \\
m &=0,1,2, \ldots, M-1 \\
N &=L M
\end{aligned}
\end{equation}

and
\begin{equation}
\begin{aligned}
\operatorname{ALIAS}_{L, m}(x) & \triangleq \sum_{l=0}^{L-1} x(m+l M) \\
m &=0,1,2, \ldots, M-1 \\
N &=L M
\end{aligned}
\end{equation}
There is a downsampling theorem introduced,
\begin{theorem}[Downsampling theorem]
    For all \(x \in  \mathbf{C}^{N}\),
    \[
        \operatorname{DOWNSAMPLE}_{L}\left( x \right) \leftrightarrow \frac{1}{L} \operatorname{ALLAS}_{L}\left( X \right)
    .\]
\end{theorem}

From this theorem, we know that downsampling in the time domain will cause the amplitudes in the frequency domain overlap, therefore extra frequencies will pollute the expected frequency bands.

    Thus before calculating the response, we need to eradicate aliasing before downsampling. We will consider to use Butterworth filters to keep the low frequency part but removing the high frequency part whose frequency is greater than the half of Nyquist frequency.

    \noindent \textbf{Butterworth filters}
    The class of continuous-time Butterworth filters\cite{SignalSystem} is that for which the magnitude squared of the frequency response \(B\left( \omega \right)\) is of the form
    \begin{equation}
        \left|B\left( \omega \right)\right|^2 = \frac{1}{1+\left( \omega / \omega_{c} \right)^{2N}}
    \label{resps}
    \end{equation}
    The parameter \(N\) is refered to as the filter order, The higher the order, the sharper the transition from passband to the stopband. The \(\omega_{c}\) is at which \(\left|B\left( \omega \right)\right|\) is at \(1 / \sqrt{2} \) its values at \(\omega = 0\).

    The magnitude squared of the frequency response of the class of discrete-time Butterworth filters is of the form
    \begin{equation}
        \left|B\left( \omega \right)\right|^2 = \frac{1}{1+ \left( \frac{\tan\left( \omega / 2 \right)}{\tan\left( \omega_{c} / 2 \right)} \right)^{2N}}
    \label{dresps}
    \end{equation}
    The frequency response \(B\left( \omega \right)\) satisfies the form
    \begin{equation}
        B\left( s \right) = \frac{\omega_{c}^{N}}{\prod_{p=1}^{N}  \left( s + s_{p} \right)}
    \label{bs}
    \end{equation}
    where
    \[
        s_{p} = \omega_{c} \exp \left\{ j \left[ \frac{\pi\left( 2p + 1 \right)}{2N} \frac{\pi}{2} \right]  \right\}
    \]
    is the pole of \(B\left( j\omega \right)B^{*}\left( j\omega \right)\).

\begin{centering}
\begin{figure}
     \centering
     \begin{subfigure}[b]{.8\textwidth}
         \centering
         \includegraphics[width=\textwidth]{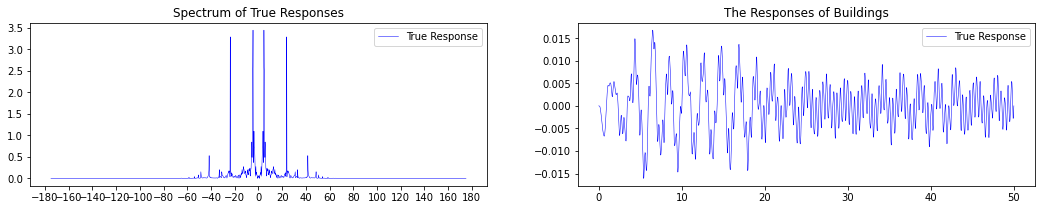}
         \caption{}
         \label{fig:smallmag}
     \end{subfigure}
     \begin{subfigure}[b]{.8\textwidth}
         \centering
         \includegraphics[width=\textwidth]{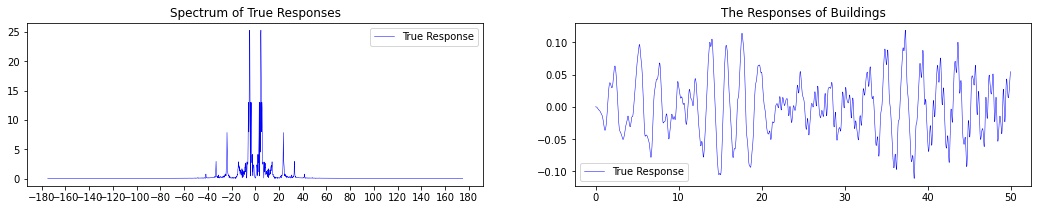}
         \caption{}
         \label{fig:midmag}
     \end{subfigure}
     \begin{subfigure}[b]{.8\textwidth}
         \centering
         \includegraphics[width=\textwidth]{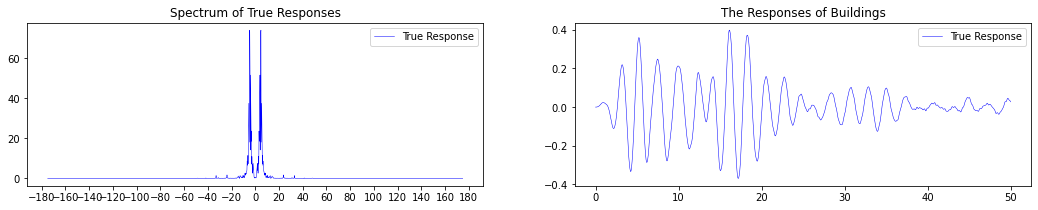}
         \caption{}
         \label{fig:largemag}
     \end{subfigure}
      \caption{Several cases that have evident discrepancy in the magnitudes. The small magnitude response has higher frequency components.}
      \label{fig:mag}
\end{figure}
\end{centering}

\begin{figure}[htpb]
        \begin{center}
            \includegraphics[width=0.5\textwidth]{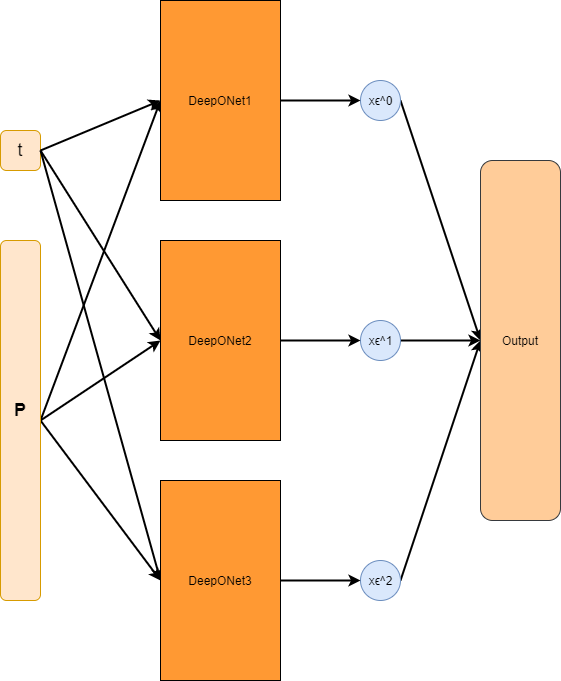}
        \end{center}
        \caption{The multiscale DeepONet}
        \label{fig:msdeeponet}
        \end{figure}

\noindent \textbf{Amplitude separation in data}
If we plot the responses, we could conclude that the high frequency has smaller magnitude, as shown in the Figure \ref{fig:mag}.
The structure now is
\begin{equation}
    f\left( P,t \right) = \sum_{i=1}^{n} \epsilon^{i}f_{i}\left( P,t \right)
\label{msdo}
\end{equation}
where \(f_{i}\)'s are DeepONet as shown in the schematics, but the trunk net of each are multiscale nets with different scales, the schematics as shown in Figure \ref{fig:msdeeponet}. Based on our observation, the data with larger magnitude contains lower frequency, the data with smaller magnitude contains higher frequency. Thus, we should configure the \(f_{i}\) with more subnets as  \(i\) increases.

        We will compare the multiscale MSDeepONet with a large MSDeepONet in section \ref{sec:Results}, which shows huge improvements.

\section{Numerical Results}
\label{sec:Results}
    \subsection{Pre-experiements}
    \label{subsec:preresults}
\begin{centering}
\begin{figure}
         \centering
         \includegraphics[width=0.75\textwidth]{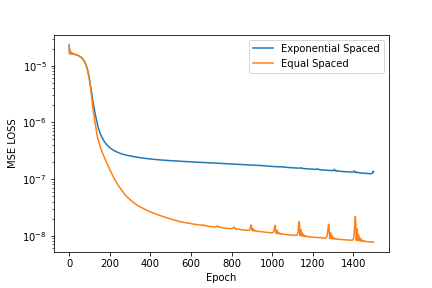}
         \caption{The training loss of neural networks with different set up of scales.}
         \label{fig:spaced}
\end{figure}
\end{centering}

    \noindent{\bf Selection of scales - \(l*n\) or  \(l^{n}\)}
    The setting of scales of subnet of multiscale deep neural network is still an open problem. We will compare two different settings of scales of the subnets of multiscale deep neural network in this subsection, one is that we set the scales as equal spaced, \( \{ l, 2l, 3l, \ldots, \kappa_{up}\}  \), the other is that we set the scale as exponential spaced \(\{s^{0}, s^{1},s^{2},\ldots, \kappa_{up} \} \), where \(l = \frac{\kappa_{up}}{N}, s = \kappa_{up}^{\frac{1}{N}}\) and \(N\) is the number of subnets. We set the branch net as fully connected neural network but the trunk net as multiscale deep neural network with different scale settings, whose capability will be proved in the next pre-experiments.  In this experiment,  a random case will be selected and an upper bound frequency \(\kappa_{up}\) will be considered, which should be the upper bound of frequency of the corresponding random case. The number of subnets is \(30\), \(\kappa_{up} = 60\times 2\pi\). As shown in Figure \ref{fig:spaced}, The loss of the equal spaced set up decreases faster and will be smaller compared with the exponentially spaced case. In addition, this case give us another bonus observation that to learn the oscillating case to a satisfactory accuracy, we need the MSE loss to be \(O\left( 10^{-6} \right)\) at least.

    \noindent{\bf Determining DeepONet structure}
    Like the previous experiment, we use only one seismic record and the corresponding response for the DeepONet structure determination.

    The multiscale neural network considered as trunk net has \(100\) four-layer subnetworks with \(10\) neurons in each layer, whose scales are \(\{1, 1+2\pi,\ldots ,1+2k\pi,\) \(  \ldots , 1+200\pi\} \) respectively, and the activation is \(\sin\left( x \right)\). The output of trunk net is 1000 dimension. The fully connected network considered as trunk net is also a four-layer neural network, using \(\sin\left( x \right)\) as activation function but with \(1000\) hidden neurons to keep the total number of hidden neurons to be the same as the multiscale neural network.

    The multiscale neural network considered as branch net has \(100\) four-layer subnetworks with \(5\) neurons in each layer, whose scales are \(\{1, 1+2\pi,\ldots ,1+2k\pi,\) \( \ldots , 1+200\pi\} \) respectively, as the same as the scales of multiscale net considered as trunk net. The activation is also \(\sin\left( x \right)\). The dimension of the output of branch net is as the same as the dimension of output of trunk net. The fully connected network considered as trunk net now shoule be a four-layer neural network with \(500\) hidden neurons whose activation function is also \(\sin\left( x \right)\).

\begin{centering}
    \begin{figure}
     \centering
     \begin{subfigure}[b]{\textwidth}
         \includegraphics[width=\hsize]{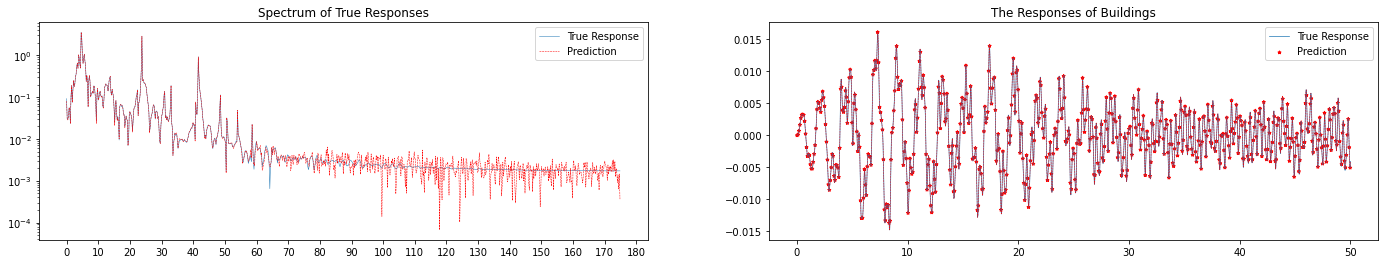}
         \caption{(bFCN-tMS)}
         \label{fig:TMS_BFCN}
     \end{subfigure}
     \hfill
     \begin{subfigure}[b]{\textwidth}
         \includegraphics[width=\hsize]{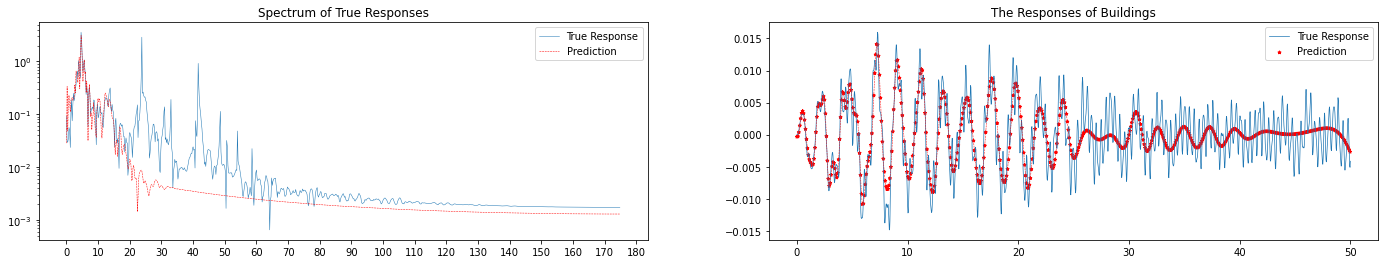}
         \caption{(bFCN-tFCN)}
         \label{fig:TFCN_BFCN}
     \end{subfigure}
     \hfill
     \begin{subfigure}[b]{\textwidth}
         \includegraphics[width=\hsize]{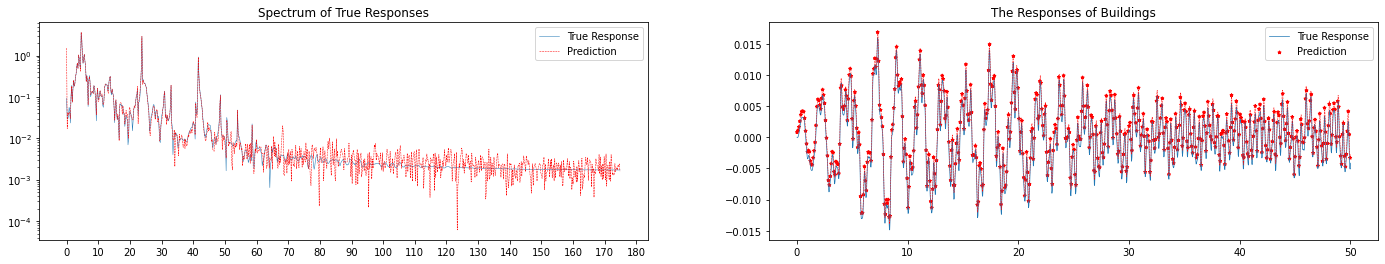}
         \caption{(bMS-tMS)}
         \label{fig:TMS_BMS}
     \end{subfigure}
      \begin{subfigure}[b]{\textwidth}
         \includegraphics[width=\hsize]{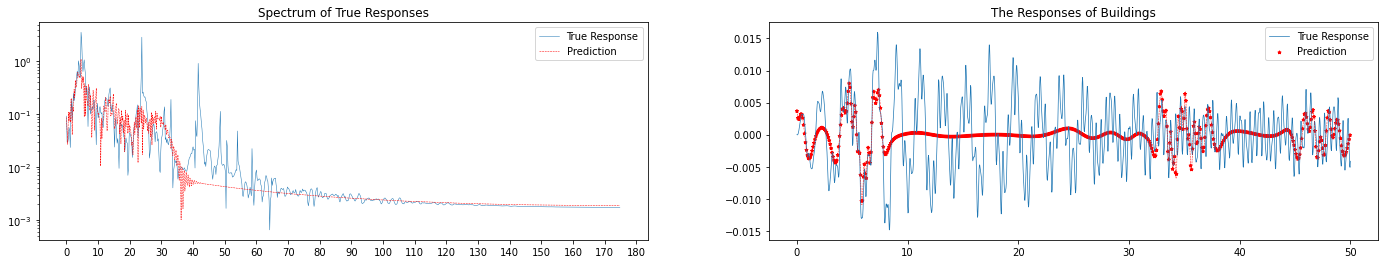}
         \caption{(bMS-tFCN)}
         \label{fig:TFCN_BFCN}
     \end{subfigure}
     \caption{The output of the training data for different network structures. As we can see, the bFCN-tMS (a) DeepONet captures the frequency the best, but the one without multiscale DNN in the trunk of DeepONet bFCN-tFCN (b) and bMS-tFCN (d) fails to do so in 100 epochs.}
        \label{fig:pre-exp}
\end{figure}
\end{centering}

    Those four DeepONets with different structures will be trained up to \(1500\) epochs. The training results can be shown in the Figure \ref{fig:pre-exp}. As we can see, both the bFCN-tMS DeepONet and the bMS-tMS DeepONet capture the frequencies, but other 2 structures only learn some low-frequency waves after 1500-epoch training. These results prove that there is no high frequency part that will affect the learning for the inputs of branch net \(\left[ P\left( t_1 \right),P\left( t_2 \right),\ldots P\left( t_{m} \right) \right] \), or at least, we would not gain any improvements if we use the scales \(\{1, 1+2\pi,\ldots ,1+2k\pi,\) \( \ldots , 1+200\pi\} \) for these inputs.

This results should be expected since solution \(x\left( t \right)\) will contain high frequency parts with respect to \(t\) if the right hand side contains high frequency part with respect to \(t\) from the Equation \(\left( \ref{Solid} \right)\). Thus the trunk net whose input is \(t\) should use multiscale deep neural network to accelerate to capture the high frequency part. Both multiscale neural net and fully connected neural net could be considered as the proper choice for branch net, but considering training the multiscale neural network will be more time-consuming, we select the fully connected neural net as branch net. We will use the bFCN-tMS structure to learn. The diagram of the bFCN-tMS structure can be shown in Fig. \ref{fig:onet_MS}.

\subsection{Results of multiscale DeepONet with amplitude separation}
\label{subsec: msdon_result}
\noindent{\bf Predict the responses for a specific floor}
The multiscale DeepONet in this section uses 3 sub-DeepONet, as shown in the Figure \(\ref{fig:msdeeponet}\). For each sub-DeepONet, the branch net is a 4-layer fully connected neural network with 128 hidden neurons in each layer, whose activation function is Relu. The trunk net is a multiscale neural network, but with different scales for different order of \(\varepsilon\). The DeepONet multiplying \(\varepsilon^2\) at the output should contain more scales since the responses whose amplitude is small contain higher frequency parts. The DeepONet multiplying \(\varepsilon^{0}\) at the output do not need contain such scales since the responses whose amplitude is large only contain lower frequency parts. We choose the scale \(\{1,1+2\pi, \ldots, 1+2k \pi, \ldots , 1+200\pi\} \) for the DeepONet multiplying \(\varepsilon^2\), the scale \(\{1,1+2\pi,\ldots ,1+2k\pi,\ldots , 1+100\pi\} \) for the DeepONet multiplying \(\varepsilon\), the scale \(\{1,1+2\pi,\ldots , 1+2k\pi,\ldots ,1+20\pi\} \) for the DeepONet multiplying \(\varepsilon^{0}\). Each subnets of those multiscale deep neural networks are 4-layer neural networks with 8 hidden neurons in each layer. As comparison, we use a large multiscale DeepONet without amplitude separation trained by the same datasets and tested by the same testing datasets. The trunk net of multiscale DeepONet without amplitude separation is a multiscale neural network which contains 100 four-layer subnets with \(24\) hidden neurons in each layer, whose scale is the same as the multiscale DeepONet with amplitude seperation. The branch net of the multiscale DeepONet without amplitude separation is the 4-layer fully connected neural network with \(3\times 128\) hidden neurons in each layer, whose activation function is also relu. The learning rate is \(10^{-3}\) and total learning procedure contains 1500 epochs with 40 batches per epoch.
\begin{centering}
\begin{figure}
     \centering
     \begin{subfigure}[b]{.8\textwidth}
         \centering
         \includegraphics[width=\textwidth]{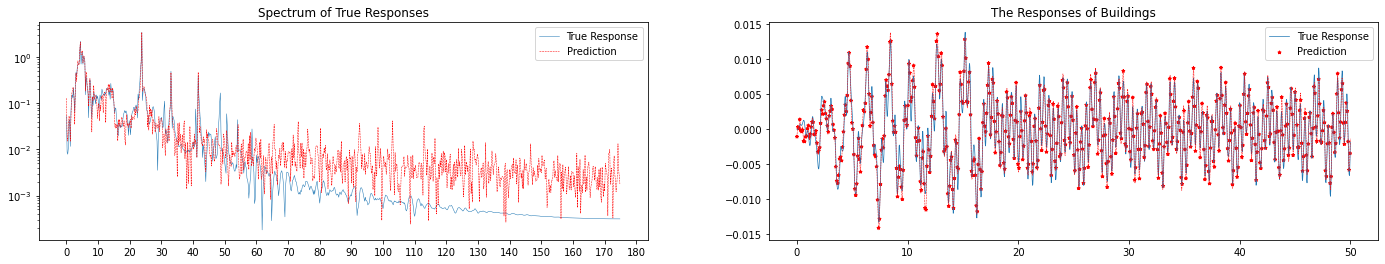}
         \caption{}
         \label{fig:smallpred}
     \end{subfigure}
     \begin{subfigure}[b]{.8\textwidth}
         \centering
         \includegraphics[width=\textwidth]{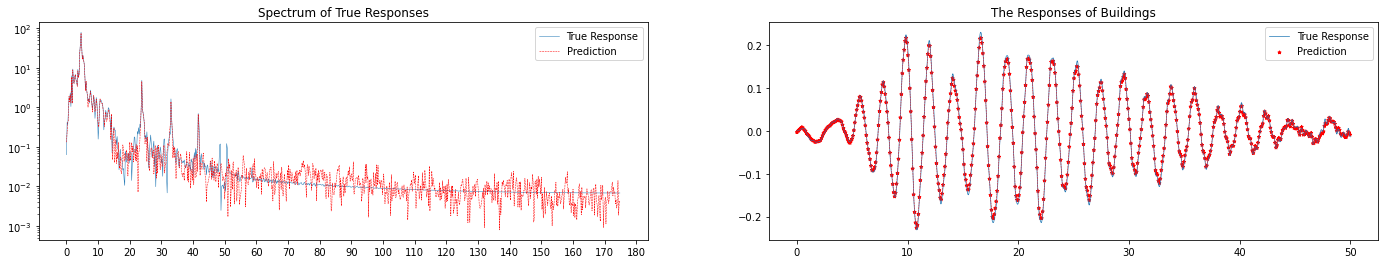}
         \caption{}
         \label{fig:largepred}
     \end{subfigure}
     \caption{The prediction of the multiscale deepONet. As shown in the figure, the neural network capture the correct frequencies.}
      \label{fig:pred}
\end{figure}
\end{centering}

\begin{centering}
\begin{figure}
     \centering
     \begin{subfigure}[b]{.45\textwidth}
         \centering
         \includegraphics[width=\textwidth]{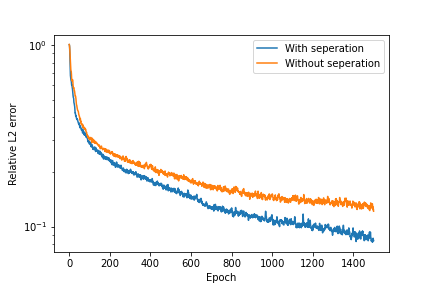}
         \caption{The relative L2 training error of the multiscale DeepONet with and without amplitude separation}
         \label{fig:loss_namps}
     \end{subfigure}
     \begin{subfigure}[b]{.45\textwidth}
         \centering
         \includegraphics[width=\textwidth]{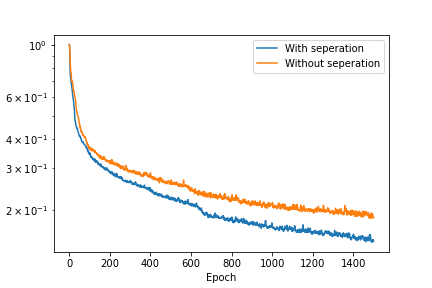}
     \caption{The relative L2 testing error of the multiscale DeepONet with and without amplitude separation}
     \label{fig:loss_nof}
     \end{subfigure}
     \caption{The multiscale deepONet with amplitude separation gives better convergence results.}
      \label{fig:loss}
\end{figure}
\end{centering}

The results of the multi scale DeepONet are shown in the Figure \ref{fig:pred}, the neural network predict the correct frequencies.
By using the data augmentation, there is no over-fitting issue and the prediction to the test cases show satisfactory accuracy. The mean relative L2 error for testing cases at last is \(0.13\).
The comparison of evolution of relative L2 error of training dataset and testing dataset for the multiscale DeepONet with amplitude separation comparing with what for multiscale DeepONet without separation are shown in Figure \ref{fig:loss}. We could conclude that the amplitude separation idea indeed have some contributions even though the data augmentation is also applied during training procedure of multiscale DeepONet without amplitude separation.

\noindent{\bf Prediction of the responses in multiple floors}
The preceding results show the multiscale DeepONet with amplitude separation could predict the response of a specific floor with a satisfactory accuracy. In this section, we show that the multiscale DeepONet with amplitude separation not only could predict the response of specific floor, but also could predict the responses for many floors simultaneously, excited by the same seismic waves. The main difference of the multiple-floor cases is that there are multiple outputs with different kinds of frequencies for different outputs. Specifically, the responses of lower floor contains higher frequencies with smaller amplitudes while the responses of higher floor contains lower frequencies with larger amplitudes. Such a amplitudes-frequency disparity can be handled with the proposed amplitude separation in the multiscale DeepONet.

\begin{centering}
\begin{figure}
    \centering
     \begin{subfigure}[b]{.8\textwidth}
         \centering
         \includegraphics[width=\textwidth]{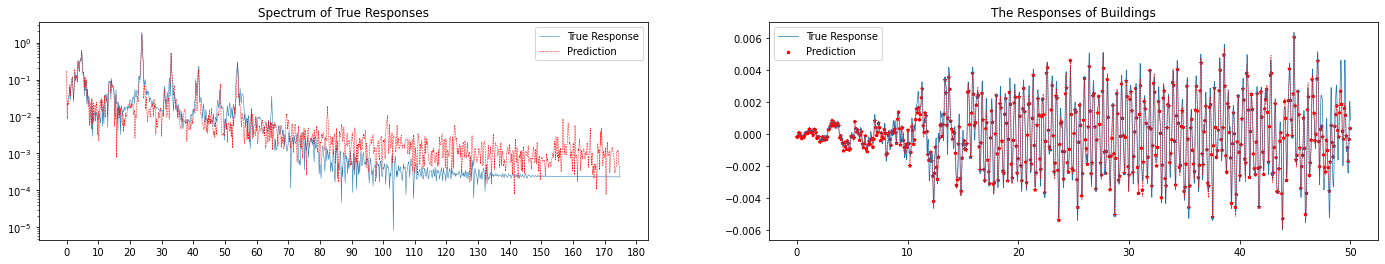}
         \caption{Floor 2}
         \label{fig:pred_F1}
     \end{subfigure}
     \begin{subfigure}[b]{.8\textwidth}
         \centering
         \includegraphics[width=\textwidth]{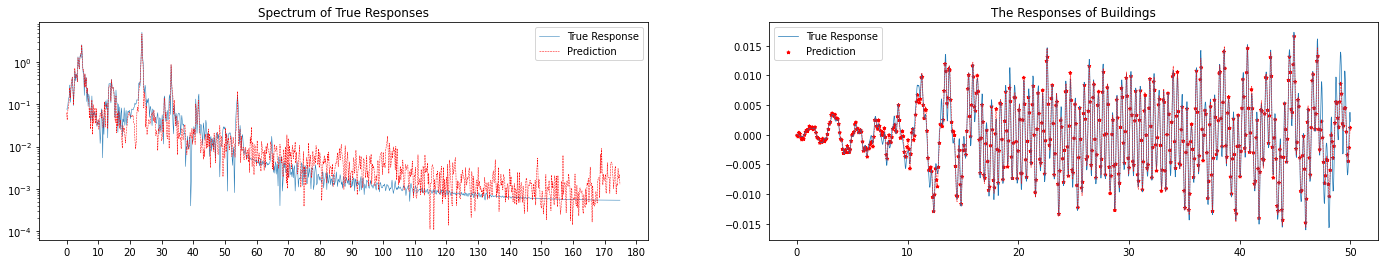}
         \caption{Floor 3}
         \label{fig:pred_F2}
     \end{subfigure}
      \begin{subfigure}[b]{.8\textwidth}
         \centering
         \includegraphics[width=\textwidth]{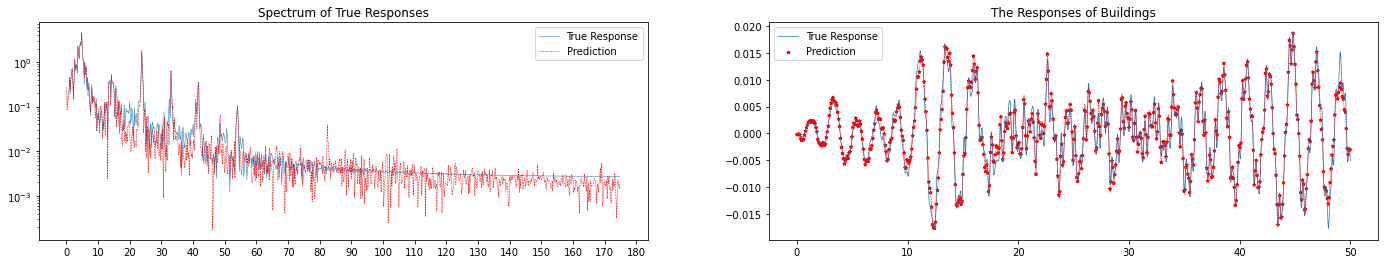}
         \caption{Floor 4}
         \label{fig:pred_F3}
     \end{subfigure}
     \begin{subfigure}[b]{.8\textwidth}
         \centering
         \includegraphics[width=\textwidth]{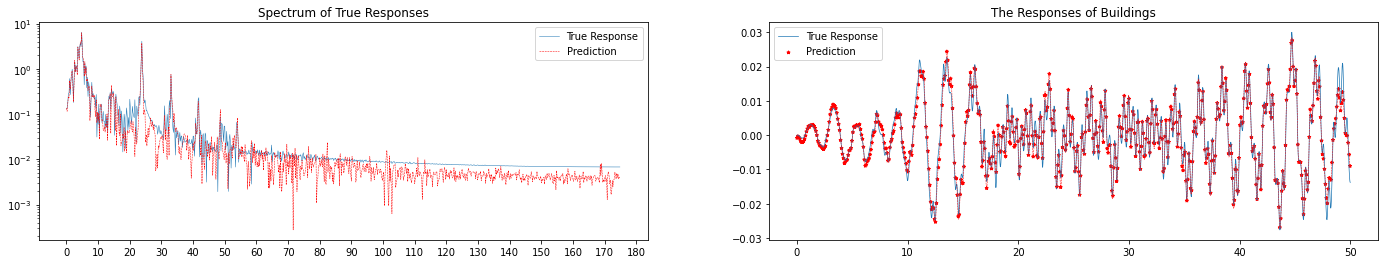}
         \caption{Floor 5}
         \label{fig:pred_F4}
     \end{subfigure}
      \begin{subfigure}[b]{.8\textwidth}
         \centering
         \includegraphics[width=\textwidth]{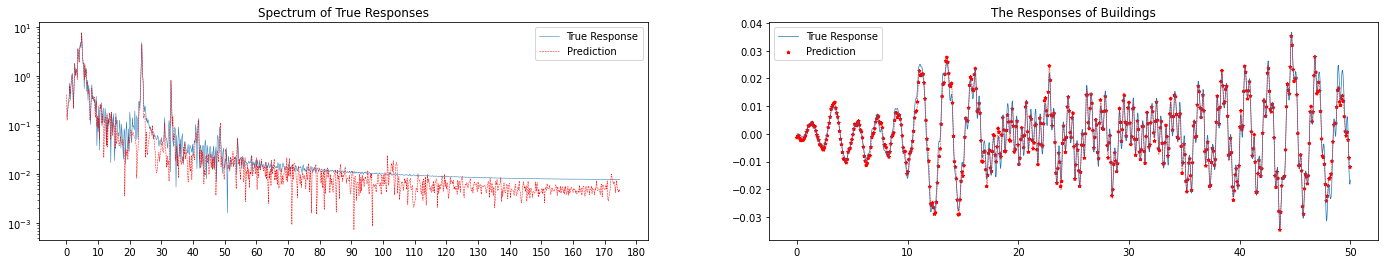}
         \caption{Floor 6}
         \label{fig:pred_F5}
     \end{subfigure}
     \begin{subfigure}[b]{.8\textwidth}
         \centering
         \includegraphics[width=\textwidth]{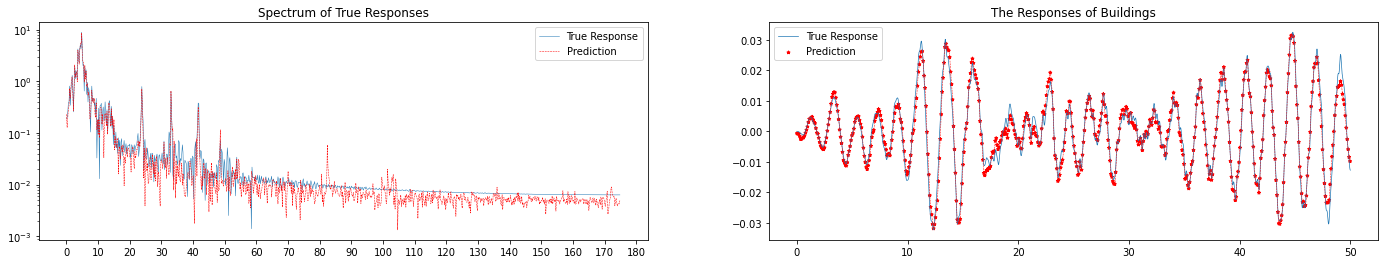}
         \caption{Floor 7}
         \label{fig:pred_F6}
     \end{subfigure}
     \begin{subfigure}[b]{.8\textwidth}
         \centering
         \includegraphics[width=\textwidth]{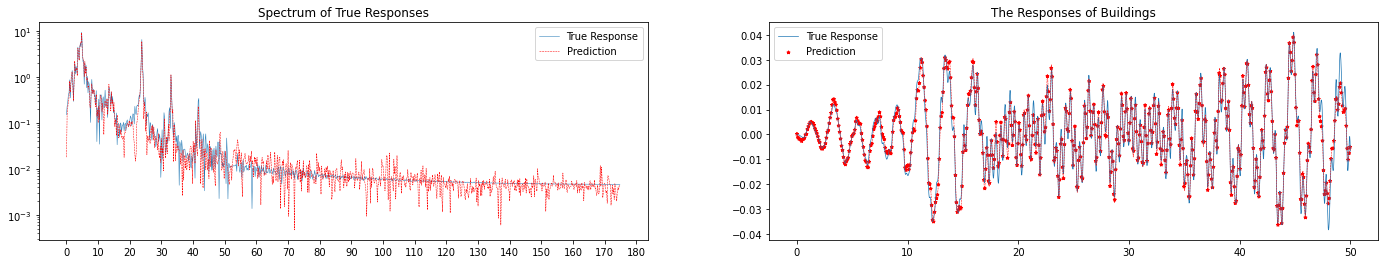}
         \caption{Floor 8}
         \label{fig:pred_F7}
     \end{subfigure}
     \caption{The predictions of the multiscale deepONet for seismic responses of multiple floors. The neural network capture the correct frequencies (left column) and time responses (right column).}
      \label{fig:pred_multi_f}
\end{figure}
\end{centering}

Now, the difference between the multiscale DeepONet here is that the output of the neural network will be an \(l\)-D array. Other settings are the same as what in the single-floor case. The results are shown in Fig. \ref{fig:pred_multi_f}. We could conclude that the multiscale DeepONet captured the desired frequencies with satisfactory accuracy from the left figures of Fig. \ref{fig:pred_multi_f} and time responses from the right figures.

\section{Conclusion}
\label{sec:conclusion}
In this paper, we proposed a multiscale DeepONet to handle the mappings between oscillating inputs and outputs. For situations with the amplitude separation in signals, we also proposed a specific structure to handle large disparity of scales in amplitudes, resulting in satisfactory results. The proposed multiscale DeepONet also shows its power when predicting the responses  containing different frequencies for different floors.

\noindent{\bf Future work} The fast convergence of multiscale neural network needs rigorous mathematical analysis even with lots of evidences of the advantages of multiscale neural networks. Meanwhile, the theory behind amplitude separation also needs further studies. The response operator \(\mathcal{R}\) function acts like a lowpass filter and is linear. For nonlinear operators such as the forward operator in the inverse medium scattering problem,  more research will be done to see the performance of the multiscale DeepONet, a success for these cases will be very helpful for solving inverse scattering problems.
\section{Acknowledgement}
The authors like to thank Prof. GE. Karniadakis for bringing the attention of this research project to our attention and Dr. Kamaljyoti Nath for helpful discussions and assistance during this work.
\bibliographystyle{unsrt}

\begin{thebibliography}{1}

\bibitem{weinan18} W. N. E,  B. Yu, The deep Ritz method: A deep learning-based numerical algorithm
for solving variational problems. Communications in Mathematics and Statistics, 6(1):1\-12, 2018.



\bibitem {han18}J. Han, A. Jentzen, W.N. E, Solving high-dimensional partial
differential equations using deep learning. Proceedings of the National
Academy of Sciences. 2018 Aug 21;115(34):8505-10.

\bibitem{gk20}  X. Jin, S. Cai, H. Li, GE. Karniadakis, NSFnets (Navier-Stokes flow nets): Physics-informed neural networks for the incompressible Navier-Stokes equations. Journal of Computational Physics. 2021 Feb 1;426:109951.

\bibitem{gk19} M. Raissi, P. Perdikaris, and G. E. Karniadakis, Physics-informed neural networks:
A deep learning framework for solving forward and inverse problems involving nonlinear
partial differential equations. Journal of Computational Physics, 378:686–707, 2019.

\bibitem{cai19} W. Cai, X.G. Li, and L.Z. Liu. A phase shift deep neural network for high frequency
approximation and wave problems. SIAM Journal on Scientific Computing. 2020;42(5):A3285-312.

\bibitem{cai20} Z.Q. Liu, W. Cai, and Z.Q. John Xu, Multi-scale Deep Neural Network (MscaleDNN) for
Solving Poisson-Boltzmann Equation in Complex Domains, arXiv:2007.11207, 2020, Communications in Computational Physics. 2020 Jun;28(5):1970-2001.

\bibitem{luDeepONet2020}
Lu~Lu, Pengzhan Jin, and George~Em Karniadakis.
\newblock {{DeepONet}}: {{Learning}} nonlinear operators for identifying
  differential equations based on the universal approximation theorem of
  operators. arXiv preprint arXiv:1910.03193. 2019 Oct 8.

\bibitem{tianpingchen1995}
{Tianping Chen} and {Hong Chen}.
\newblock Universal approximation to nonlinear operators by neural networks with arbitrary activation functions and its application to dynamical systems. IEEE Transactions on Neural Networks. 1995 Jul;6(4):911-7.

\bibitem{OpenSystemEarthquake}
Open {{System}} for {{Earthquake Engineering Simulation}} - {{Home Page}}.


\bibitem{NumericalSimulationLargescale2006}
Jun-Hong Ding, Xian-Long Jin, Yi-Zhi Guo, and Gen-Guo Li.
\newblock Numerical simulation for large-scale seismic response analysis of
  immersed tunnel.
\newblock Engineering structures. 2006 Aug 1;28(10):1367-77.

\bibitem{chopraDynamicsStructures2011}
Anil~K. Chopra.
\newblock {\em Dynamics of {{Structures}}}.
\newblock {Pearson}, 4th edition edition.


\bibitem{li2021} Z.Y. Li, N. Kovachki, K. Azizzadenesheli, B. Liu, K. Bhattacharya, A. Stuart, A. Anandkumar, Fourier neural operator for parametric partial differential equations, arXiv preprint arXiv:2010.08895, 2021.

\bibitem{WaveletbasedSimulation2017}
V.~L. Nithin, S.~Das, and H.~B. Kaushik.
\newblock Wavelet-based simulation of scenario-specific nonstationary
  accelerograms and their {{GMPE}} compatibility.
\newblock Soil Dynamics and Earthquake Engineering. 2017 Aug 1;99:56-67.

\bibitem{MDFT}
Julius~O. Smith.
\newblock {\em Mathematics of the Discrete Fourier Transform ({{DFT}})}.
\newblock {http://ccrma.stanford.edu/ jos/mdft/}, accessed (date accessed).

\bibitem{SignalSystem}
Alan~V. Oppenheim, Alan~S. Willsky, and S.~Hamid Nawab.
\newblock {\em Signals \& Systems (2nd Ed.)}.
\newblock {Prentice-Hall, Inc.}

\end{thebibliography}

\end{document}